\documentclass[11pt,reqno,a4paper,twoside]{amsart}
\usepackage{bbm,amsmath,amsfonts,amssymb,latexsym,color}
\definecolor{gr1}{gray}{0.85}
\definecolor{gr2}{gray}{0.5}
\definecolor{gr3}{gray}{0.7}
\definecolor{pink}{rgb}{1,0.6,0.6}
\definecolor{lgr}{rgb}{0.6,1,0.6}
\definecolor{lbl}{rgb}{0.6,0.6,1}
\definecolor{lye}{rgb}{1,1,0.6}

\newtheorem{proposition}{Proposition}[section]
\newtheorem{corollary}[proposition]{Corollary}

\newtheorem{remark}[proposition]{Remark}
\newtheorem{remarks}[proposition]{Remarks}

\newcounter{z}
\newenvironment{enum}{\setcounter{z}{0}
\begin{list}{{\rm\Alph{z}.}}{\usecounter{z}
\setlength{\topsep}{1ex}\setlength{\labelwidth}{0.5cm}
\setlength{\leftmargin}{0.65cm}\setlength{\labelsep}{0.25cm}
\setlength{\parsep}{0pt}}}{\end{list}}
\newcounter{z1}

\newcounter{z2}
\newenvironment{enum2}{\setcounter{z2}{0}
\begin{list}{\arabic{z2})}{\usecounter{z2}
\setlength{\topsep}{0ex}\setlength{\labelwidth}{0.5cm}
\setlength{\leftmargin}{0.65cm}\setlength{\labelsep}{0.25cm}
\setlength{\parsep}{0pt}}}{\end{list}}

\newcommand{\bea}{\begin{eqnarray*}}
\newcommand{\eea}{\end{eqnarray*}}

\newcommand{\brm}{\begin{rm}}
\newcommand{\erm}{\end{rm}}

\newcommand{\sn}{\mathcal{S}_n}

\newcommand{\dd}{\mathsf{r}}
\newcommand{\uu}{\mathsf{u}}
\newcommand{\des}{\mathsf{des}}
\newcommand{\inv}{\mathsf{inv}}
\newcommand{\lrmax}{\mathsf{lrmax}}
\newcommand{\maxdes}{\mathsf{maxdes}}
\newcommand{\maj}{\mathsf{maj}}
\newcommand{\lbsum}{\mathsf{lbsum}}
\newcommand{\ba}{\mbox{-}}
\newcommand{\std}[1]{{#1}^{\ast}}

\begin{document}

\title{The area above the Dyck path of a permutation}
\author{Mark Dukes, Astrid Reifegerste}
%\begin{keyword} Permutation \sep Dyck path \sep  Simion-Schmidt \sep Bruhat poset 
%\MSC 05A05, 05A17, 05A19
%\end{keyword}
\address{Science Institute, University of Iceland, Reykjav\'{i}k, Iceland}
\email{dukes@hi.is}
%\ead{dukes@hi.is}
\address{Faculty of Mathematics, University of Magdeburg, Germany}
\email{astrid.reifegerste@ovgu.de}
%\ead{astrid.reifegerste@ovgu.de}

\maketitle

\begin{abstract}
In this paper we study a mapping from permutations to Dyck paths. 
A Dyck path gives rise to a (Young) diagram and we give relationships between statistics on permutations and statistics on their corresponding diagrams. The distribution of the size of this diagram is discussed and a generalization given of a parity result due to Simion and Schmidt. We propose a filling of the diagram which determines the permutation uniquely. Diagram containment on a restricted class of permutations is shown to be related to the strong Bruhat poset.
\end{abstract}

\section{Statistics on permutations and their paths}

Let $\sn$ be the set of all permutations $\pi=\pi_1\pi_2\ldots\pi_n$ on $[n]=\{1,2,\ldots,n\}$. 
Let $B=\{b_1,\ldots , b_n\}$ with $b_1< \cdots < b_n$ be any finite subset of $\mathbb{N}$.
Define the {\em{standardization}} of a permutation $\pi$ on $B$ to be the permutation $\std{\pi} \in \sn$ 
obtained from $\pi$ by replacing $b_i$ with the integer $i$. For example $\std{3\,5\,9\,4}=1\,3\,4\,2$.
In this paper we will use the language of generalized patterns. 
The reader unfamiliar with generalized patterns should consult Claesson~\cite{claesson}.

We will define a mapping $D$ that maps permutations to Dyck paths.
This mapping can be traced back to Knuth \cite[\S 2.2.1 Exercise 3, 4 and 5]{knu}
where it arose as a result of of analysing push and pop operations in a stack and the admissible sequences of these operations.
Using more recent terminology, the mapping mentioned in \cite{knu} is from the collection of $(2\ba 3\ba 1)$-avoiding permutations 
to the collection of Dyck words.

Define the recursive map $D$ from $\sn$ to the set of Dyck paths of length $2n$ as follows:
Set $D(\emptyset)=\emptyset$ and $D(1)=\uu\dd$. 
Given a permutation $\pi\in\sn$, split $\pi$ into $\pi=\pi_Ln\pi_R$ and define $D(\pi)=\uu D(\std{\pi_L})\dd D(\std{\pi_R})$. 

For example, the permutation $\pi=5\,3\,1\,4\,8\,2\,7\,6\in\mathcal{S}_8$ is mapped to the path $D(\pi)=\uu D(4\,2\,1\,3)\dd D(1\,3\,2)=\uu D(4\,2\,1\,3)\dd\uu\uu\dd\dd\uu\dd=\uu\uu\dd\uu\uu\dd\uu\dd\dd\dd\uu\uu\dd\dd\uu\dd$.

\begin{center}
\unitlength0.35cm
\begin{picture}(9,9)
\put(0,0){\color{gr2}\line(1,1){8}}
\put(0,1){\color{gr2}\line(1,0){1}}
\put(0,2){\color{gr2}\line(1,0){2}}
\put(0,3){\color{gr2}\line(1,0){3}}
\put(0,4){\color{gr2}\line(1,0){4}}
\put(0,5){\color{gr2}\line(1,0){5}}
\put(0,6){\color{gr2}\line(1,0){6}}
\put(0,7){\color{gr2}\line(1,0){7}}
\put(0,8){\color{gr2}\line(1,0){8}}
\put(0,0){\color{gr2}\line(0,1){8}}
\put(1,1){\color{gr2}\line(0,1){7}}
\put(2,2){\color{gr2}\line(0,1){6}}
\put(3,3){\color{gr2}\line(0,1){5}}
\put(4,4){\color{gr2}\line(0,1){4}}
\put(5,5){\color{gr2}\line(0,1){3}}
\put(6,6){\color{gr2}\line(0,1){2}}
\put(7,7){\color{gr2}\line(0,1){1}}
\linethickness{0.7pt}
\put(0,0){\line(0,1){2}}\put(0,2){\line(1,0){1}}
\put(1,2){\line(0,1){2}}\put(1,4){\line(1,0){1}}
\put(2,4){\line(0,1){1}}\put(2,5){\line(1,0){3}}
\put(5,5){\line(0,1){2}}\put(5,7){\line(1,0){2}}
\put(7,7){\line(0,1){1}}\put(7,8){\line(1,0){1}}
\put(0.4,-0.2){\makebox(0,0)[cc]{\tiny\sf 0}}
\put(1.4,0.8){\makebox(0,0)[cc]{\tiny\sf 2}}
\put(2.4,1.8){\makebox(0,0)[cc]{\tiny\sf 4}}
\put(3.4,2.8){\makebox(0,0)[cc]{\tiny\sf 6}}
\put(4.4,3.8){\makebox(0,0)[cc]{\tiny\sf 8}}
\put(5.4,4.8){\makebox(0,0)[cc]{\tiny\sf 10}}
\put(6.4,5.8){\makebox(0,0)[cc]{\tiny\sf 12}}
\put(7.4,6.8){\makebox(0,0)[cc]{\tiny\sf 14}}
\put(8.4,7.8){\makebox(0,0)[cc]{\tiny\sf 16}}
\end{picture}
\end{center}
\vspace*{2ex}

For the usual representation of a Dyck path rotate the figure clockwise by $\pi/4$.
We call a pair of entries $\dd\uu$ in a Dyck path a {\it{valley}}.

Let us note that the mapping $D$ with domain $\sn$ is not a bijection, 
as is easily seen since $D(1\, 3\, 2)=D(2\,3\,1)=\uu\uu\dd\dd\uu\uu$.
The permutation $\pi$ is also associated to a partition $\Lambda(\pi)$ 
whose Young diagram arises when lengthening the ascents and descents of the path. 
The Young diagram of any partition obtained from a permutation $\pi\in\sn$ fits in the shape $(n-1,n-2,\ldots,1)$. 
In the following, we will identify a partition with its Young diagram and vice versa.
To recover the length of $\pi$ from a given partition, 
we allow parts equal to 0 and require that $\Lambda(\pi)$ has exactly $n-1$ parts.
In the above example, we have $\Lambda(\pi)=(7,5,5,2,1,1,0)$.

It transpires that Knuth's map relates the classical statistics of partitions with natural permutation statistics.
For a permutation $\pi\in\sn$, we call the numbers $a_1,a_2,\ldots,a_n$ 
the {\it left border numbers} of $\pi$ where $a_i$ is the position of the 
rightmost element to the left of $\pi_i$ in $\pi$ which is greater than $\pi_i$, 
or 0 if no such element exists. 
For instance, the left border numbers of $\pi=5\,3\,1\,4\,8\,2\,7\,6$ are $0,1,2,1,0,5,5,7$. 
Clearly, $a_i=0$ whenever $\pi_i$ is a left-to-right maximum of $\pi$ 
(that is, an element which is greater than all elements to its left), 
and $a_i=i-1$ whenever $i-1$ is a descent of $\pi$ (that is, $\pi_{i-1}>\pi_i$). 
Furthermore, each number $a_i>0$ is necessarily a descent of $\pi$ and so the nonzero left border numbers form the descent set of $\pi$.   

\begin{proposition}
For $\pi\in\sn$, let $a$ be the sequence of the left border numbers 
$a_2,a_3,\ldots,a_n$ of $\pi$ in decreasing order. Then $\Lambda(\pi)=a$.
\end{proposition}

\begin{proof}
First we note that the valleys ($\dd\uu$) of the path $D(\pi)$ are in bijection with the descents of $\pi$. 
Let $\pi_i=n$. We may assume that $i<n$; otherwise we consider the reduced permutation $\pi_1\pi_2\ldots\pi_{n-1}$. 
(The identity permutation $12\ldots n$ is mapped to the path $\uu^n\dd^n$ having no valleys.) 
Clearly, $i$ is a descent of $\pi$. 
By construction, we have $D(\pi)=\uu D(\pi_1'\ldots\pi_{i-1}')\dd D(\pi_{i+1}'\ldots\pi_n')$ where 
$\pi_L^\ast=\pi_1'\ldots\pi_{i-1}'\in\mathcal{S}_{i-1}$ and $\pi_R^\ast=\pi_{i+1}'\ldots\pi_n'\in\mathcal{S}_{n-i}$.
Obviously, $\pi_L^\ast$ and $\pi_R^\ast$ together have exactly one descent less than $\pi$. 
On the other hand, since $D(\pi_R^\ast)$ starts with $\uu$, the paths 
$D(\pi_L^\ast)$ and $D(\pi_R^\ast)$ together have exactly one valley less than $D(\pi)$. 
Consequently, the descents of $\pi$ corresponds to the valleys of Knuth's path. 
More exactly, the vertical projection of the valley point corresponding to descent $i$ meets the diagonal at point $2i$. 
Hence the descents of $\pi$ are the different nonzero parts of $\Lambda(\pi)$.

The multiplicity of part $i>0$ in $\Lambda(\pi)$ is given by the number of 
consecutive $\uu$'s following the step $\dd$ which marks the valley associated to $i$. 
Consider $D(\pi)=\uu\ldots\dd D(\pi_{i+1}'\ldots\pi_k')\ldots$. 
The number $a_{i+1}=i$ accounts for the initial $\uu$ step of $D(\pi_{i+1}'\ldots\pi_k')$. 
Splitting $\pi_{i+1}'\ldots\pi_k'$ into $LmR$ with maximal element $m$ yields 
$D(\pi_{i+1}'\ldots\pi_k')=\uu D(L^\ast)\dd D(R^\ast)$. 
Repeating this process, for every decomposition the initial sequence of consecutive 
$\uu$'s is extended by a further $\uu$ while $L$ is nonempty. 
In any case, $\pi_i$ is the rightmost element to the left of $m$ which is greater than $m$. 
Thus the number of $j$ satisfying $a_j=i$ equals the number of $\uu$'s we identify to be the multiplicity of part $i$.        
\end{proof}

\begin{remarks} \label{codefromshape}
\brm
\hspace*{1ex}\\[-4ex]
\begin{enum}
\item Note that the left border numbers $a_2,\ldots,a_n$ can be recovered in their original order from $\Lambda(\pi)$ ($a_1=0$ by definition). Given the partition, we set $a_{i+1}=i$ for each nonzero part $i$ of $\Lambda(\pi)$. The remaining parts of $\Lambda(\pi)$ are then placed as follows. Replace the undetermined $a_j$'s (from left to right) by the greatest possible part which is smaller than $j$.

For the above example, we have $\Lambda(\pi)=(7,5,5,2,1,1,0)$, and hence $a_2=1$, $a_3=2$, $a_6=5$, and $a_8=7$. The number $a_4$ equals to the greatest remaining part of $\Lambda(\pi)$ which is smaller than $4$, namely 1. Then $a_5$ takes the greatest remaining part which is smaller than $5$, namely 0, and lastly $a_7$ needs to be 5.
\item Analogously, we can define the {\it right border number} $b_i$ as the smallest integer $j>i$ with $\pi_i<\pi_j$, or $n+1$ if there is no such element. The numbers $n+1-b_n,n+1-b_{n+1},\ldots,n+1-b_1$ are just the left border numbers of the reverse of $\pi$.
\end{enum}
\erm
\end{remarks}

Given a partition $\lambda$ as in the diagram below,
we define the {\em{size of a reversed hook}} of a cell in $\lambda$ to be 1 plus the number of cells above it, plus the number of cells to its left.
Let the {\em{maximal reverse hook length}} be the reverse hook length that is maximal over all cells of the partition.
In the diagram below there are two cells having maximal reverse hook length which is 7: 
the lowest and rightmost pink cell, and the rightmost blue cell.

For a fixed partition $\lambda$, the number of permutations $\pi\in\sn$ with $\Lambda(\pi)=\lambda$ can be given in terms of the length of particular reverse hooks. 
For this, we divide $\lambda$ successively into disjoint rectangles as follows: 
Consider the corner with maximal reverse hook length. 
If there is more than one, choose the leftmost. 
Outline the rectangle consisting of all cells left and above the corner 
and repeat it for 
the subshape below this rectangle and the subshape to the right of this rectangle.
The process is carried out until all cells have an assigned rectangle. 
(Finally, we have as many rectangles as $\lambda$ has corners.) 
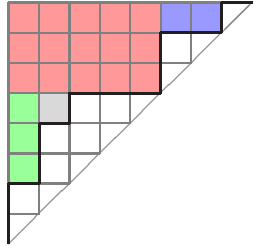
\begin{center}
\unitlength0.4cm
\begin{picture}(8,9)
\put(0,5){\fcolorbox{pink}{pink}{\makebox(5,3){}}}
\put(1,4){\fcolorbox{gr1}{gr1}{\makebox(1,1){}}}
\put(0,2){\fcolorbox{lgr}{lgr}{\makebox(1,3){}}}
\put(5,7){\fcolorbox{lbl}{lbl}{\makebox(2,1){}}}
\put(0,0){\color{gr2}\line(1,1){8}}
\put(0,1){\color{gr2}\line(1,0){1}}
\put(0,2){\color{gr2}\line(1,0){2}}
\put(0,3){\color{gr2}\line(1,0){3}}
\put(0,4){\color{gr2}\line(1,0){4}}
\put(0,5){\color{gr2}\line(1,0){5}}
\put(0,6){\color{gr2}\line(1,0){6}}
\put(0,7){\color{gr2}\line(1,0){7}}
\put(0,8){\color{gr2}\line(1,0){8}}
\put(0,0){\color{gr2}\line(0,1){8}}
\put(1,1){\color{gr2}\line(0,1){7}}
\put(2,2){\color{gr2}\line(0,1){6}}
\put(3,3){\color{gr2}\line(0,1){5}}
\put(4,4){\color{gr2}\line(0,1){4}}
\put(5,5){\color{gr2}\line(0,1){3}}
\put(6,6){\color{gr2}\line(0,1){2}}
\put(7,7){\color{gr2}\line(0,1){1}}
\linethickness{0.7pt}
\put(0,0){\line(0,1){2}}\put(0,2){\line(1,0){1}}
\put(1,2){\line(0,1){2}}\put(1,4){\line(1,0){1}}
\put(2,4){\line(0,1){1}}\put(2,5){\line(1,0){3}}
\put(5,5){\line(0,1){2}}\put(5,7){\line(1,0){2}}
\put(7,7){\line(0,1){1}}\put(7,8){\line(1,0){1}}
\end{picture}
\end{center}
Let $w_i$ and $h_i$ be the width and height, respectively, of the rectangle whose bottom right-hand corner is contained in the $i$th column of $\lambda$. 

\begin{proposition} \label{partition_into_rectangle}
For a partition $\lambda\subseteq(n-1,n-2,\ldots,1)$ with corners in the columns $i_1,i_2\ldots,i_k$ there are 
$$
\prod_{j=1}^k {w_{i_j}+h_{i_j}-1\choose w_{i_j}-1}
$$
permutations $\pi\in\sn$ with $\Lambda(\pi)=\lambda$.
\end{proposition}

\begin{proof}
As mentioned above, the corner in column $i$ corresponds to the descent $i$. By construction, $w_i-1$ is the number of elements between the rightmost larger element to the left of $\pi_i$ (which is on position $a_i$) and $\pi_i$. But the height $h_i$ counts the elements between $\pi_i$ and the leftmost larger element to the right of $\pi_i$ (which is on position $b_i$). If we know all the elements arising somewhere between $a_i$ and $b_i$ for each descent $i$ and furthermore, if we know which of these elements arise to the left of position $i$, then we can determine the permutation $\pi$ completely. 
\end{proof}

\begin{remark}
\brm
Knuth's path construction is closely related to the representation of permutations as decreasing binary trees (see \cite[pp. 23]{sta} for a slightly modified version). Given a permutation $\pi=\pi_Ln\pi_R\in\sn$, define $T(\pi)$ to be the tree with root $n$ having left and right subtrees $T(\pi_L^*)$ and $T(\pi_R^*)$, respectively, obtained by removing $n$. This yields inductively a tree such that the left predecessor of a vertex $\pi_i$ is the rightmost element to the left of $\pi_i$ which is greater than $\pi_i$, that is, $\pi(a_i)$, while the right predecessor of $\pi_i$ is the leftmost element greater than $\pi_i$ to the right of $\pi_i$, that is, $\pi(b_i)$. 
\erm
\end{remark}

Some of the natural statistics of $\Lambda(\pi)$ can be expressed in terms of statistics of $\pi$. Let $\des(\pi)$ denote the number of descents of $\pi$, $\lrmax(\pi)$ the number of left-to-right maxima of $\pi$, $\maj(\pi)$ the major index of $\pi$ (that is, the sum of all elements in the descent set $\{i:\pi_i>\pi_{i+1}\}$), and $\lbsum(\pi)$ the sum of the left border  numbers of $\pi$.

\begin{corollary} \label{statistics}
Let $\pi\in\sn$ and $\Lambda$ its associated partition. Then 
\begin{enumerate}
\item[(i)] The number of different nonzero parts of $\Lambda$ is equal to $\des(\pi)$. 
\item[(ii)] The number of nonzero parts of $\Lambda$ is equal to $n-\lrmax(\pi)$.  
\item[(iii)] The largest part of $\Lambda$ is equal to the largest descent of $\pi$.
\item[(iv)] The sum of all different parts of $\Lambda$ is equal to $\maj(\pi)$.   
\item[(v)] The sum of all parts (equivalently, area) of $\Lambda$ is equal to $\lbsum(\pi)$.  
\end{enumerate}
\end{corollary}

By the definition of the left border numbers, $\lbsum(\pi)$ is the sum of the initial position of the occurrences of all the (consecutive) patterns $(k+1)\omega k$ with $\omega\in\mathcal{S}_{k-1}$ arbitrary. But we also can express this statistic in terms of the number of occurrences of certain patterns (without considering the positions where the occurrences appear in $\pi$).

For any pattern $\sigma$, let $(\sigma)\pi$ be the number of the occurrences of $\sigma$ in $\pi$. 
Furthermore, we use $1\ba\bar{4}\ba3\ba\bar{4}\ba2$ to denote a pattern $1\ba3\ba2$ 
which is part of neither a $1\ba4\ba3\ba2$ nor a $1\ba3\ba4\ba2$ pattern. 

\begin{proposition} \label{cp-pattern}
We have $\lbsum(\pi)=(2\ba1)\pi+(1\ba\bar{4}\ba3\ba\bar{4}\ba2)\pi$ for any permutation $\pi$.
\end{proposition}

\begin{proof}     
For any $j$, the number $a_j$ counts the number of integers $i<j$ for which there exists $k$ such that $i\le k<j$ and $\pi_k>\pi_j$. (This is compatible with the original definition of $a_j$. Suppose $a_j=m$, that is, $m$ is maximal with $m<j$ and $\pi_m>\pi_j$. Then for $i=1,\ldots,m$ there exists such a $k$, namely $k=m$, whereas for $m+1,\ldots,j-1$ there is no such $k$.)

In the case of $\pi_i>\pi_j$, the number of inversions $(2\ba1)\pi$ will count the relevant pairs $(i,j)$. If $\pi_i<\pi_j$ then there must exist a $k$ such that $\pi_i\pi_k\pi_j$ is an occurrence of the pattern $1\ba3\ba2$. To avoid the multiple counting of $(i,j)$, we assume that $k$ is maximal, that is, we must rule out elements larger than $\pi_k$ in the subword $\pi_i\pi_{i+1}\ldots\pi_{j-1}\pi_j$. Hence the number of such pairs $(i,j)$ is $(1\ba\bar{4}\ba3\ba\bar{4}\ba2)\pi$.      
\end{proof}

Since the second expression in Proposition \ref{cp-pattern} counts particular occurrences of the pattern $1\ba3\ba2$, the statistic $\lbsum(\pi)$ coincides with the number of inversions of $\pi$ if $\pi$ avoids $1\ba3\ba2$. In this case, $\Lambda(\pi)$ is just the conjugated permutation diagram, see \cite{rei}. Thus Knuth's map induces a bijection between $1\ba3\ba2$-avoiding permutations and Dyck paths which takes the inversion number to the area statistic which was already observed by Bandlow and Killpatrick \cite{bankil}.

\section{Distribution of lbsum}

In \cite[Prop.~1]{simschm}, Simion and Schmidt show that there are as many even as odd $1\ba3\ba2$-avoiding permutations in $\sn$ if $n$ is even while the excess of even over odd $1\ba3\ba2$-avoiding permutations equals a Catalan number if $n$ is odd.
We give the following generalization.

\begin{proposition} \label{zig:zag}
For $n$ even, the number of permutations in $\sn$ for which $\lbsum$ takes an even value is the same as the number of permutations in $\sn$ for which the value is odd. For $n$ odd, the numbers differ by the tangent numbers.
\end{proposition}

\begin{proof}
By definition, for a permutation $\pi=\pi_Ln\pi_R\in\sn$ the path $D(\pi)$ is the concatenation of $\uu D(\pi_L^\ast)\dd$ and $D(\pi_R^\ast)$. 
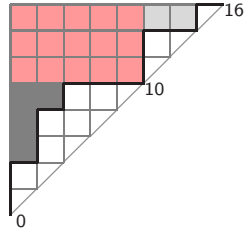
\begin{center}
\unitlength0.35cm
\begin{picture}(9,8.5)
\put(0,5){\fcolorbox{pink}{pink}{\makebox(5,3){}}}
\put(5,7){\fcolorbox{gr1}{gr1}{\makebox(2,1){}}}
\put(0,2){\fcolorbox{gr2}{gr2}{\makebox(1,3){}}}
\put(1,4){\fcolorbox{gr2}{gr2}{\makebox(1,1){}}}
\put(0,0){\color{gr2}\line(1,1){8}}
\put(0,1){\color{gr2}\line(1,0){1}}
\put(0,2){\color{gr2}\line(1,0){2}}
\put(0,3){\color{gr2}\line(1,0){3}}
\put(0,4){\color{gr2}\line(1,0){4}}
\put(0,5){\color{gr2}\line(1,0){5}}
\put(0,6){\color{gr2}\line(1,0){6}}
\put(0,7){\color{gr2}\line(1,0){7}}
\put(0,8){\color{gr2}\line(1,0){8}}
\put(0,0){\color{gr2}\line(0,1){8}}
\put(1,1){\color{gr2}\line(0,1){7}}
\put(2,2){\color{gr2}\line(0,1){6}}
\put(3,3){\color{gr2}\line(0,1){5}}
\put(4,4){\color{gr2}\line(0,1){4}}
\put(5,5){\color{gr2}\line(0,1){3}}
\put(6,6){\color{gr2}\line(0,1){2}}
\put(7,7){\color{gr2}\line(0,1){1}}
\linethickness{0.7pt}
\put(0,0){\line(0,1){2}}\put(0,2){\line(1,0){1}}
\put(1,2){\line(0,1){2}}\put(1,4){\line(1,0){1}}
\put(2,4){\line(0,1){1}}\put(2,5){\line(1,0){3}}
\put(5,5){\line(0,1){2}}\put(5,7){\line(1,0){2}}
\put(7,7){\line(0,1){1}}\put(7,8){\line(1,0){1}}
\put(0.4,-0.2){\makebox(0,0)[cc]{\tiny\sf 0}}
\put(5.4,4.8){\makebox(0,0)[cc]{\tiny\sf 10}}
\put(8.4,7.8){\makebox(0,0)[cc]{\tiny\sf 16}}
\end{picture}
\end{center}
Interpreting $\lbsum(\pi)$ as area of the shape $\Lambda(\pi)$, it is easy to see that $\lbsum(\pi)=\lbsum(\pi_L^\ast)+\lbsum(\pi_R^\ast)+k(n-k)$ where $\pi_k=n$. (Note that the path $D(\pi)$ first returns to the diagonal at $2k$.)

In the case of even $n$, the map $\phi:\pi\mapsto\pi_Rn\pi_L$ is an involution on $\sn$ which switches the parity of $\lbsum$ since $n$ is placed on position $n+1-k$ in $\phi(\pi)$.

If $n$ is odd, $\lbsum(\pi)$ is an even number if and only if $\lbsum(\pi_L^\ast)$ and $\lbsum(\pi_R\ast)$ have the same parity. Let $e_n$ (resp. $o_n$) be the number of permutations in $\sn$ with an even (resp. odd) left border number sum. Then we have
\bea
e_n=\sum_{k=1}^n{n-1\choose k-1}\big(e_{k-1}e_{n-k}+o_{k-1}o_{n-k}\big)\\
o_n=\sum_{k=1}^n{n-1\choose k-1}\big(e_{k-1}o_{n-k}+o_{k-1}e_{n-k}\big)
\eea
(where $e_0=1,\;o_0=0$). Because $e_k=o_k$ for each even $k$, we obtain
$$
\Delta_n:=e_n-o_n=\sum_{k=2\atop\mbox{\tiny $k$ even}}^{n-1}{n-1\choose k-1}\Delta_{k-1}\Delta_{n-k}\qquad(\Delta_1=1).
$$
\end{proof}

The considerations made at the beginning of the previous proof immediately yield a recursion formula for the ordinary generating function of the statistic $\lbsum$. 

\begin{proposition} 
The function $F_n(x)=\sum\limits_{\pi\in\sn} x^{\lbsum(\pi)}$ satisfies
\begin{equation} \label{recursion}
F_n(x)=\sum_{k=1}^n{n-1\choose k-1}F_{k-1}(x)F_{n-k}(x)x^{k(n-k)}
\end{equation}
where $F_0(x)=1$.
\end{proposition}

As a consequence, we obtain the following information on the expectation of the random 
variable $X_n:\sn\to\Bbb{N}_0$ with $X_n(\pi)=\lbsum(\pi)$ where the probability measure is the uniform distribution on $\sn$. We also have a representation for the variance in terms of harmonic numbers. 

\begin{corollary}
For $n\ge 2$, we have 
\bea
{\mathsf E}(X_n)&=&(n+1)(\mbox{$\frac{n}{2}$}-H_{n,1})+n\\
Var(X_n)&=&2n(n+2)-(n+1)H_{n,1}-(n+1)^2H_{n,2}
\eea
where $H_{n,m}=1+\frac{1}{2^m}+\ldots+\frac{1}{n^m}$ denotes the harmonic number of order $m$. 
\end{corollary}
%\newpage

\begin{remarks}
\brm
\hspace*{1ex}\\[-4ex]
\begin{enum}
\item As mentioned in Section 1, the statistic $\lbsum$ is identical with the inversion number $\inv$ for $1\ba3\ba2$-avoiding permutations. Over $\sn(1\ba3\ba2)$, the generating function of $\inv$ is the $q$-Catalan polynomial $C_n(q)$ which satisfies the recursion 
$$
C_n(q)=\sum_{k=1}^n C_{k-1}(q)C_{n-k}(q)q^{k-1}
$$ 
as shown by F\"urlinger and Hofbauer \cite{furhof}. Equation (\ref{recursion}) generalizes this result.
\item The random variable $X_n$ has the same distribution function as the random variable ${n\choose 2}-Y_n$ where $Y_n$ has appeared in literature several times, see \cite{knuth_maa}, \cite{kirsch}, \cite{hwang}. It is known to measure the major cost of the in-situ permutation algorithm or the left path length of random binary search trees. Kirschenhofer et al. \cite{kirsch} and Hwang and Neininger \cite{hwang} have given limit laws and general moments for  $Y_n$ (and therefore $X_n)$.
\end{enum}
\erm
\end{remarks}

Refining, let $G_n(x,y,p,q)$ be the generating function of the statistic quadruple $(\lbsum,\des,\maxdes,\lrmax)$ where $\maxdes(\pi)$ is defined to be the largest descent of $\pi$. Note that, for a given permutation $\pi\in\sn$, these statistics measure the area, the number of corners, the number of columns, and $n$ minus the number of rows of the associated shape $\Lambda(\pi)$. Using the decomposition of $\pi$ into $\pi_Ln\pi_R$, we obtain for all $n\ge2$  
\bea
G_n(x,y,p,q)&=&G_{n-1}(x,y,p,q)\\
&&+\sum_{k=1}^{n-1} {n-1\choose k-1}G_{k-1}(x,y,1,q)G_{n-k}(x,y,p,1)x^{k(n-k)}yp^kq
\eea
where $G_0(x,y,p,q)=G_1(x,y,p,q)=1$. The generating function of these polynomials 
\bea
g(x,y,p,q,z)&=&\sum_{n\ge0}\mbox{$\frac{1}{n!}$}\,x^{n\choose 2}z^nG_n(x^{-1},y,p,q)
\eea
admits the recursion
$$
\frac{\partial g(x,y,p,q,z)}{\partial z}=g(x,y,p,q,xz)-ypg(x,y,1,q,xpz)(1-g(x,y,p,1,qz)).
$$
Particular instances of this may be solved; for example, one easily finds $g(-1,1,1,1,z)=1+\tanh(z)$, 
giving an alternative proof of Proposition \ref{zig:zag}.

\section{A permutation representation}

In this section, we provide $\Lambda(\pi)$ with an additional filling to obtain a unique assignment between permutations and filled shapes.

Given the shape $\Lambda(\pi)$, label its columns from left to right with 1 to $k$ (where $k$ is the last descent of $\pi$). Then label the rows of length $i\ge1$ from bottom to top with the numbers $j$ (in increasing order) for which $a_j=i$. (Note that the bottom row of length $i$ has label $i+1$.) Lastly, put a dot into the cell in column $i$ and row $j$ whenever $(i,j)$ is an inversion of $\pi$. For $\pi=5\,3\,1\,4\,8\,2\,7\,6$, the procedure yields 
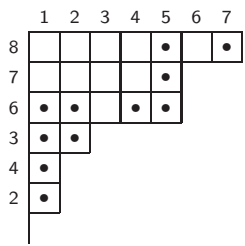
\begin{center}
\unitlength0.4cm
\begin{picture}(8,8)
\put(1,0){\line(0,1){7}}\put(2,1){\line(0,1){6}}
\put(3,3){\line(0,1){4}}\multiput(4,4)(1,0){3}{\line(0,1){3}}\multiput(7,6)(1,0){2}{\line(0,1){1}}
\multiput(1,6)(0,1){2}{\line(1,0){7}}\multiput(1,4)(0,1){2}{\line(1,0){5}}
\put(1,3){\line(1,0){2}}\multiput(1,1)(0,1){2}{\line(1,0){1}}
\put(1.5,7.5){\makebox(0,0)[cc]{\tiny\sf 1}}
\put(2.5,7.5){\makebox(0,0)[cc]{\tiny\sf 2}}
\put(3.5,7.5){\makebox(0,0)[cc]{\tiny\sf 3}}
\put(4.5,7.5){\makebox(0,0)[cc]{\tiny\sf 4}}
\put(5.5,7.5){\makebox(0,0)[cc]{\tiny\sf 5}}
\put(6.5,7.5){\makebox(0,0)[cc]{\tiny\sf 6}}
\put(7.5,7.5){\makebox(0,0)[cc]{\tiny\sf 7}}
\put(0.5,6.5){\makebox(0,0)[cc]{\tiny\sf 8}}
\put(0.5,5.5){\makebox(0,0)[cc]{\tiny\sf 7}}
\put(0.5,4.5){\makebox(0,0)[cc]{\tiny\sf 6}}
\put(0.5,3.5){\makebox(0,0)[cc]{\tiny\sf 3}}
\put(0.5,2.5){\makebox(0,0)[cc]{\tiny\sf 4}}
\put(0.5,1.5){\makebox(0,0)[cc]{\tiny\sf 2}}
\put(1.5,1.5){\circle*{0.3}}\put(1.5,2.5){\circle*{0.3}}\put(1.5,3.5){\circle*{0.3}}\put(1.5,4.5){\circle*{0.3}}
\put(2.5,3.5){\circle*{0.3}}\put(2.5,4.5){\circle*{0.3}}\put(4.5,4.5){\circle*{0.3}}\put(5.5,4.5){\circle*{0.3}}
\put(5.5,5.5){\circle*{0.3}}\put(5.5,6.5){\circle*{0.3}}\put(7.5,6.5){\circle*{0.3}}
\end{picture}
\end{center}
Recovering the permutation from the filled shape is an easy matter: By counting the dots in the $i$th column, we know the number of inversions $(i,\cdot)$ which is enough to determine $\pi_i$ successively.

To describe the tableaux resulting from a permutation in this way, we first focus on the minimal filling. Consider the splitting of the shape introduced in Section 1. It is easy to see that each cell which arises in the rightmost column of some rectangle must be filled. As explained in the proof of Proposition \ref{partition_into_rectangle}, these cells correspond to the elements to the right of a descent top $\pi_i$ before the first occurrence of an element which is greater than $\pi_i$.  
 
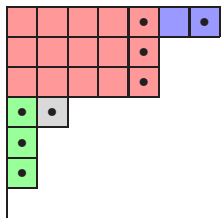
\begin{center}
\unitlength0.4cm
\begin{picture}(8,8)
\put(0,4){\fcolorbox{pink}{pink}{\makebox(5,3){}}}
\put(1,3){\fcolorbox{gr1}{gr1}{\makebox(1,1){}}}
\put(0,1){\fcolorbox{lgr}{lgr}{\makebox(1,3){}}}
\put(5,6){\fcolorbox{lbl}{lbl}{\makebox(2,1){}}}
\put(0,0){\line(0,1){7}}\put(1,1){\line(0,1){6}}
\put(2,3){\line(0,1){4}}\multiput(3,4)(1,0){3}{\line(0,1){3}}\multiput(6,6)(1,0){2}{\line(0,1){1}}
\multiput(0,6)(0,1){2}{\line(1,0){7}}\multiput(0,4)(0,1){2}{\line(1,0){5}}
\put(0,3){\line(1,0){2}}\multiput(0,1)(0,1){2}{\line(1,0){1}}
\multiput(0.5,1.5)(0,1){3}{\circle*{0.3}}
\put(1.5,3.5){\circle*{0.3}}
\multiput(4.5,4.5)(0,1){3}{\circle*{0.3}}
\put(6.5,6.5){\circle*{0.3}}
\end{picture}
\end{center}

But these are exactly those cells which have to be filled absolutely. The unique permutation in bijection with this tableau has the property that for all $i$, if $j>i$ is the smallest integer with $\pi_i<\pi_j$ then $\pi_i<\pi_k$ for all $k\ge j$. This means that $\pi$ avoids the pattern $2\ba3\ba1$.

As mentioned in Section 1, the other extreme case (all the cells are filled) corresponds to $1\ba3\ba2$-avoiding permutations which proves the following result.

\begin{proposition}
There is a bijection between $1\ba3\ba2$-avoiding and $2\ba3\ba1$-avoiding permutations in $\sn$ which preserves the corresponding shape and hence all the statistics mentioned in Corollary \ref{statistics}.
\end{proposition}

By the tableau construction, it is an easy matter to read the number of occurrences of $1\ba3\ba2$ from the filling. If we have an empty cell and to its right a filled cell, precisely the empty cell in column $i$ and the filled cell in column $j$, within the row labeled with $k$ then the sequence $\pi_i\pi_j\pi_k$ is an occurrence of $1\ba3\ba2$ in $\pi$. (Note that $i<j<k$ by the definition of the labels, and $\pi_i<\pi_k<\pi_j$ by the definition of the cell filling.) Conversely, if $\pi_i\pi_j\pi_k$ is an $1\ba3\ba2$-occurrence then $a_k\ge j$ and the row labeled with $k$ comprises at least $j$ cells. Therefore the tableau records all the occurrences of $1\ba3\ba2$ in $\pi$.

\begin{remark}
\brm
The occurrences of $2\ba3\ba1$ are not immediately countable. Two filled cells (in columns $i$ and $j$) in the same row $k$ mean that $\pi_i\pi_j\pi_k$ is either an occurrence of $2\ba3\ba1$ or $3\ba2\ba1$. To find out which case is true one has to check whether the tableau contains the cell $(i,j)$ (in column $i$ and row $j$). If so, then the cell is filled and $\pi_i\pi_j\pi_k$ is a decreasing subsequence; otherwise we have an occurrence of $2\ba3\ba1$.    
\erm
\end{remark}

Furthermore, it turns out that the poset of $1\ba3\ba2$-avoiding permutations ordered by containment of the associated partitions is just the (strong) Bruhat poset restricted to $\sn(1\ba3\ba2)$.

\begin{proposition} \label{poset}
For $\pi,\sigma\in\sn(1\ba3\ba2)$, we have $\Lambda(\pi)\subset\Lambda(\sigma)$ if and only if $\pi<\sigma$ in the Bruhat order.  
\end{proposition}

\begin{proof}     
Let $\pi$ be covered by $\sigma$ in the Bruhat poset, that is,
$$
\sigma=\pi_1\ldots\pi_{i-1}\pi_j\pi_{i+1}\ldots\pi_{j-1}\pi_i\pi_{j+1}\ldots\pi_n
$$
with $\pi_i<\pi_j$ such that $\sigma$ has exactly one inversion more than $\pi$. Consequently, there is no $k$ with $i<k<j$ and $\pi_i<\pi_k<\pi_j$. Because $\pi$ avoids $1\ba3\ba2$, we even have $\pi_k<\pi_i$ for all $k=i+1,\ldots,j-1$. The figure shows the regions (shaded) in which the dots $(k,\pi_k)$ have to be placed to guarantee that no $1\ba3\ba2$ occurs in $\pi$ or $\sigma$. 
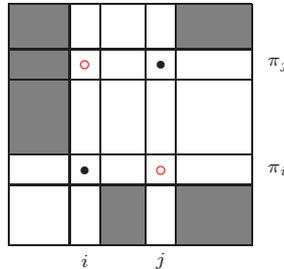
\begin{center}
\unitlength0.4cm
\begin{picture}(8,8)
\put(0,3){\fcolorbox{gr2}{gr2}{\makebox(2,5){}}}
\put(3,0){\fcolorbox{gr2}{gr2}{\makebox(1.5,2){}}}
\put(5.5,0){\fcolorbox{gr2}{gr2}{\makebox(2.5,2){}}}
\put(5.5,6.5){\fcolorbox{gr2}{gr2}{\makebox(2.5,1.5){}}}
\multiput(0,0)(8,0){2}{\line(0,1){8}}
\multiput(0,0)(0,8){2}{\line(1,0){8}}
\multiput(0,2)(0,1){2}{\line(1,0){8}}
\multiput(4.5,0)(1,0){2}{\line(0,1){8}}
\multiput(0,5.5)(0,1){2}{\line(1,0){8}}
\multiput(2,0)(1,0){2}{\line(0,1){8}}
\put(2.5,2.5){\circle*{0.3}}\put(5,6){\circle*{0.3}}
\put(2.5,6){\color{red}{\circle{0.3}}}\put(5,2.5){\color{red}{\circle{0.3}}}
\put(2.5,-0.5){\makebox(0,0)[cc]{\tiny$i$}}
\put(5,-0.5){\makebox(0,0)[cc]{\tiny$j$}}
\put(8.5,6){\makebox(0,0)[lc]{\tiny$\pi_j$}}
\put(8.5,2.5){\makebox(0,0)[lc]{\tiny$\pi_i$}}
\end{picture}
\end{center}
Now we compare the left border numbers $a_1,\ldots,a_n$ of $\pi$ 
(which are the parts of $\Lambda(\pi)$) with the left border numbers $a'_1,\ldots,a'_n$ of $\sigma$ (which are the parts of $\Lambda(\sigma)$). Obviously, we have $a'_i=a_j$, $a_i=i-1$ and $a'_j=i$. All the other left border numbers do not change when exchanging $\pi_i$ and $\pi_j$, that is, $a_k=a'_k$ for $k\not=i,j$. This is trivial for $k<i$. For $k=i+1,\ldots,j-1$, if $a_k=i$ (otherwise $a_k=l$ for some $l>i$), then $a'_k=i$ since $\pi_j>\pi_i$. Finally, for $k=j+1,\ldots,n$ we have $\pi_k>\pi_j$ or $\pi_k<\pi_i$. In the first case, $a_k=a'_k=l$ for some $l<i$. If we have $a_k=j$ in the second case (the alternative is $a_k=l$ for some $l>j$), then $a'_k=j$ as well since $\pi_k<\pi_i$.
Consequently, $\Lambda(\sigma)$ is obtained from $\Lambda(\pi)$ by adding one (corner) cell.

Conversely, assume now that $\Lambda(\pi)\subset\Lambda(\sigma)$ and the both shapes differ by exactly one cell. By Proposition \ref{cp-pattern}, the number of cells in the $i$th column of the shape is equal to the number of inversions $(i,j)$ with $i<j$ in the corresponding permutation. Hence $\sigma$ is obtained from $\pi$ by a singleton reduction.     
\end{proof} 

\begin{remark}
\brm
In general the previous statement fails. The shapes of Bruhat ordered permutations do not have to be contained in each other (for example, $\Lambda(1243)=(3,0,0)$ and $\Lambda(1423)=(2,2,0)$). Conversely, permutations whose shapes are contained in each other are not necessarily comparable in Bruhat order (for example, $\Lambda(1342)=(3,0,0)$ and $\Lambda(2143)=(3,1,0)$).
\erm
\end{remark}

{\sc Open questions}

\begin{enum2}
\item If the shape of one permutation is contained inside the shape of another, then what aspect of the permutations does this reflect?
\item What are the rules for filling tableaux as introduced in Section 3?
\end{enum2}

\section*{Acknowledgments}
The authors would like to thank Einar Steingr\'imsson for planting the seed that led to this paper
and Alois Panholzer for valuable comments and pointers to the literature.

\end{document}